\documentclass[11pt]{article}   
\usepackage{amsthm, amsmath, amssymb}   
\theoremstyle{plain}

\newtheorem{theorem}{Theorem}[section]
\newtheorem{lemma}[theorem]{Lemma}
\newtheorem{proposition}[theorem]{Proposition}
\newtheorem{corollary}[theorem]{Corollary}

\theoremstyle{definition}

\newtheorem{definition}[theorem]{Definition}
\newtheorem{example}[theorem]{Example}

\theoremstyle{remark}

\newtheorem{remark}[theorem]{Remark}

\def\qed{\hfill  \framebox(5,5){}}
\def\para{\vspace{1mm}}

\def\card{{\rm Card}}

\def\cc{{\cal C}}
\def\conch{\mathfrak{C}(\cc,A,d)}

\def\conchsimple{\mathfrak{C}(\cc)}

\def\bbsimple{\mathfrak{B}(\cc)}
\def\gp{\mathfrak{G}(\cc,A{\cal P})}
\def\K{\Bbb  K}

\def\cP{{\cal P}}

\def\card{{\rm Card}}

\def\cc{{\cal C}}

\def\gp{\mathfrak{G}(\cP)}

\def\conch{\mathfrak{C}(\cc,A,d)}

\def\conchsimple{\mathfrak{C}(\cc)}

\def\bbsimple{\mathfrak{B}(\cc)}

\def\K{\Bbb  K}

\def\rdf{{\sc rdf}\,}

\title{Rational Conchoids  of Algebraic Curves.\footnote{Both authors supported by the Spanish `` Ministerio de
Educaci\'on e Innovaci\'on" under the Project MTM2008-04699-C03-01}}
\author{J.R.  Sendra \\
Dpto. de Matem\'aticas \\
        Universidad de Alcal\'a \\
        E-28871 Madrid, Spain \\rafael.sendra@uah.es
        \and J.  Sendra \\
        Dpto. Matem\'atica Aplicada  a la\\
        I.T. de Telecomunicaci\'on. \\
        E.U.I.T. Telecomunicaci\'on \\
        Universidad
Polit\'ecnica de
 Madrid, Spain \\jsendra@euitt.upm.es
 }

\date{}

\begin{document}
\maketitle

\begin{abstract}
We study the rationality of the components of the  conchoid  to  an
irreducible algebraic affine plane
 curve, excluding the trivial cases of the isotropic lines, of the lines
 through the focus and the circle centered at the focus and radius the distance involved in the conchoid.
 We prove that conchoids having all their components rational
can only be generated by rational curves. Moreover, we show that
reducible conchoids to  rational curves have always their two
components rational. In addition, we prove that the rationality of
the conchoid component, to a rational curve, does depend on the base
curve and on the focus but not on the distance.  Also, we provide an
algorithm that analyzes the rationality of all the components of the
conchoid and, in the affirmative case, parametrizes them.  The
algorithm only uses a proper parametrization of the base curve and
the focus and, hence, does not require the previous computation of
the conchoid. As a corollary, we show that the conchoid to the
irreducible conics, with conchoid-focus on the conic, are rational
and we give parametrizations. In particular we parametrize the {\it
Lima\c{c}ons of Pascal}. We also parametrize the conchoids of {\it
Nicomedes}. Finally, we show to find the focuses from where the
conchoid  is rational or with two rational components.
\end{abstract}

\vspace*{-5 mm}

\section{Introducci\'on}
The conchoid  is a classical geometric construction. Intuitively
speaking, if $\cc$ is a plane curve  (the base curve), $A$ a fixed
point in the plane (the focus), and $d$ a non-zero fixed field
element (the distance), the conchoid of $\cal C$ from the focus $A$
at distance $d$ is the (closure of) set of points $Q$ in the line
$AP$ at distance $d$ of a point $P$ varying in the curve $\cc$.
 The two classical and most famous conchoids are the Conchoid of {\it Nicomedes} ($\cc$ is a line and $A\not\in \cc$) and the {\it Lima\c{c}cons of Pascal} ($\cc$ is a circle and $A\in \cc$).
Conchoids are useful in many applications as construction of
buildings, astronomy, electromagnetic research, physics, optics,
engineering in medicine and biology, mechanical in fluid processing,
etc (see the introduction of \cite{SS08} for references).

\para

In this paper, we deal with the problem of analyzing the rationality
of the components of a conchoid and, in the affirmative case, the
actual computation of parametrizations. Clearly,  the problem can be
approached by computing the implicit equation of the conchoid to
afterwards factor it to finally apply to each factor any
parametrization algorithm. Nevertheless, we want to avoid all these
computations solving the problem  directly from the input base curve
and the focus. For this purpose, similarly as in \cite{SS08},  we
work over an algebraically closed field  $\K$ of characteristic
zero, and  curves are considered reduced; that is, they are the zero
set in ${\Bbb K}^2$ of non-constant square-free polynomials of
${\Bbb K}[y_1,y_2]$. Furthermore, if a curve is defined by the
square-free polynomial $f$, when we speak about its components, we
mean the curves defined by the non-constant irreducible factors
(over $\Bbb K$) of $f$ (see \cite{libro} for further details).

\para

In \cite{SS08} we presented a theoretical analysis of the concept
and main properties of conchoids to irreducible curves (see Section
\ref{sec-def-conchoid} for a brief summary). In this analysis, three
different types of curves have an exceptional  behavior: the
isotropic lines $y_1\pm \sqrt{-1}\, y_2=0$  (its conchoid is empty),
the  circle centered at the focus and radius the distance involved
in the conchoid (its conchoid has a zero-dimensional component) and
the lines through the focus (all conchoid components are special;
see  Section \ref{sec-def-conchoid} for this concept). For all the
other cases,  the most remarkable property  in \cite{SS08} is that
the  conchoid is a plane algebraic curve with at most two component,
being at least one of them simple (see  Section
\ref{sec-def-conchoid} for the notion of simple component).

\para

In this paper, we exclude w.l.o.g. the above  three exceptional
types of curves. In this situation, we prove that conchoids having
all their components rational can only be generated by rational
curves.  Moreover, we show that reducible conchoids to  rational
curves have always their two components rational; we call this case
{\sf double rationality}. Furthermore, we characterize rational
conchoids and double rational conchoids. From these results, one
deduces that   the rationality of the conchoid component, to a
rational curve, does depend on the base curve and on the focus but
not on the distance. To approach the problem we use similar ideas to
those in \cite{ASS96} introducing the notion of {\sf
reparametrization curve} (see Def. \ref{def-gp}) as well as the
notion of \rdf  parametrization (see Def. \ref{def-rpn}). The \rdf
concept allows us to detect the double rationality while the
reparametrization  curve is a much simpler curve than the conchoid,
directly computed from the input rational curve and the focus, and
that behaves equivalently as the conchoid in terms of rationality.
As a consequence of these theoretical results we provide an
algorithm to solve the problem. Given a proper parametrization of
the base curve and the focus,  the algorithm  analyzes the
rationality of all the components of the conchoid and, in the
affirmative case, parametrizes them.  We note that the algorithm
does not require the computation of the conchoid. In addition, we
show that the conchoid to the irreducible conics, with
conchoid-focus on the conic, are rational and we give
parametrizations. In particular we parametrize the {\it
Lima\c{c}ons of Pascal}. We also parametrize the conchoids of {\it
Nicomedes}. Finally, we show to find the focuses from where the
conchoid  is rational or with two rational components.

\section{Preliminaries on Conchoids and General Assumptions.}\label{sec-def-conchoid}
In this section we recall the notion of conchoid as well as its main
properties. For further details, we refer to \cite{SS08}. Let  $\K$
be an algebraically closed field of characteristic zero. In  $\K^2$
we consider the symmetric bilinear form
$$\mathfrak{b}((x_{1},x_{2}),(y_{1},y_{2}))= x_{1}y_{1}+x_{2}y_{2},$$ which
induces a metric vector space with light cone   of isotropy
$\mathfrak{L}=\{P\in \K^2\,|\, \mathfrak{b}(P,P)=0\}$ (see
\cite{ST71}). That is, $\mathfrak{L}$ is the union of the two lines
defined by $x_1\pm \sqrt{-1}\, x_2=0$. In this context, the circle
of center $P\in \K^2$ and radius $d\in \K$ is the plane curve
defined by $\mathfrak{b}(\bar{x}-P,\bar{x}-P)=d^2$, with
$\bar{x}=(x_1,x_2)$. We  say that the distance between  $P,Q \in
\K^2$ is $d\in \K$ if $P$ is on the circle of center $Q$ and radius
$d$. The notion of ``distance" is hence defined up to multiplication
by $\pm 1$. On the other hand, if $P\in {\Bbb K}^{2}$ is not
isotropic (i.e. $P\not\in \mathfrak{L}$) we denote by $\| P\|$ any
of the elements in $\K$  such that $\|P\|^2=\mathfrak{b}(P,P)$, and
if $P\in {\Bbb K}^2$ is isotropic, then $\|P\|=0$. In this paper we
usually work with both solutions of $\|P\|^2=\mathfrak{b}({P},{P})$.
For this reason we use the notation $\pm \|{P}\|$.

\para

In this situation,  let $\cc$ be the affine irreducible plane curve
 defined by the irreducible polynomial
 $f(\bar{y})\in \K[\bar{y}]$, $\bar{y}=(y_1,y_2)$, let $d\in \K^*$
be a non-zero field element, and let $A=(a,b)\in \K^2$. We consider
the {\sf (conchoid) incidence variety}
\[ \bbsimple=
\left\{ (\bar{x}, \bar{y},\lambda) \in {{\Bbb K}}^{2} \times {{\Bbb
K}}^{2}\times {\Bbb  K} \left/
\begin{array}{l}
f(\bar{y})=0   \\
\|\bar{x}-\bar{y}\|^2=d^2 \\
 \bar{x}=A+\lambda(\bar{y}-A)
 \end{array}
\right \} \right.
\]
and the {\sf incidence diagram}
\[
 \begin{array}{ccc} \\
& \bbsimple
 \subset  {\Bbb  K}^{2} \times {\Bbb  K}^{2} \times {\Bbb  K} & \\
& \hspace*{0.5 cm} \mbox{ $\begin{array}{c} \vspace*{-2 mm}
\pi_{1}  \\
\\
\end{array}$ \hspace{-2 mm} {\Huge $\swarrow$}
} \hspace{2 mm} \mbox{ {\Huge $\searrow$} $\begin{array}{c}
\vspace*{-2 mm}
\pi_{2} \\
\\
\end{array}$} &   \\
& \hspace*{-5 mm} \pi_{1}(\bbsimple)
 \subset
{\Bbb  K}^{2} \hspace{1 cm}    {\cal C}
 \subset {\Bbb  K}^{2} & \\ & &  \\
\end{array} \]
where
\[
\begin{array}{llllllllll}
\pi_{1}:&  {\Bbb  K}^{2} \times {\Bbb K}^{2}\times {\Bbb K}
&\longrightarrow {\Bbb K}^{2},& & & \pi_{2}: & {\Bbb  K}^{2} \times
{\Bbb K}^{2}\times  {\Bbb
K} &\longrightarrow {\Bbb K}^{2}\\
& (\bar{x}, \bar{y},\lambda) & \longmapsto \bar{x} & & & & (\bar{x},
\bar{y},\lambda) & \longmapsto \bar{y}.
\end{array}
\]
Then, we define the {\sf conchoid of
 ${\cal C}$ from the focus
$A$ and distance $d$} as the algebraic Zariski   closure in ${\Bbb
K}^{2}$ of $\pi_1(\bbsimple)$, and we denote it by $\conchsimple$;
i.e.
\[ \conchsimple= \overline{\pi_{1}(\bbsimple)}. \]
For details on how to compute the conchoid see \cite{SS08}. In
general, $A$ and $d$  are just precise elements in $\K^2$ and
$\K^*$, respectively. When this will not be the case (for instance
in
 Section \ref{sec-param-focus})  the conchoid will be denoted by
$\conch$  instead of $\conchsimple$ to emphasize this fact.

\para

\noindent Throughout this paper, we {\bf assume} w.l.o.g. that:
\begin{enumerate}
\item $\cc$  is
none of the two  lines defining the light cone of isotropy
$\mathfrak{L}$. This ensures that $\conchsimple\neq \emptyset$.
\item $\cc$ is not a circle centered at $A$ and radius $d$. If $\cc$ is such a circle, then $\conchsimple$ decomposes as the focus union the circle centered at $A$ and radius $2d$.
    This assumption avoids that the conchoid has zero-dimensional components  (compare to Theorem \ref{th-comp}).
\item $\cc$ is not a line through the  focus. If $\cc$ is such a line, then $\conchsimple=\cc$. This assumption avoids that the conchoid has all components special (compare to Theorem \ref{th-properties-simple-special}).
\end{enumerate}

\para

The following theorem (see Theorem 1  in \cite{SS08}) states the
main property on conchoids


\begin{theorem}\label{th-comp} $\conchsimple$  has at most two components and all of them have dimension 1.
\end{theorem}


Now, we recall the notion of simple and special components of a
conchoid that, as shown in \cite{SS08}, play an important role when
studying the rationality. More precisely, an irreducible component
${\cal M}$ of $\conchsimple$ is called {\sf simple} if there exists
a non-empty Zariski dense subset $\Omega \subset {\cal M}$ such
that, for $Q\in \Omega,$ $\card(\pi_{2}(\pi_{1}^{-1}(Q)))=1$.
Otherwise $\cal M$ is called {\sf special}. The next theorem states
the main property on  the existence of simple components (see
Theorem  3 in \cite{SS08}).


\begin{theorem}\label{th-properties-simple-special}
 $\conchsimple$ has  at least one simple component.
\end{theorem}


The next lemma (see Lemma 5 in \cite{SS08}) connects the
birationality of the maps in incidence diagram and the simple
components of the conchoid.


\begin{lemma}\label{lemma-pi-birracional}
 Let $\pi_{1}, \pi_{2}$  be the projections in the incidence
 diagram of $\cc$, and $\cal M$ an irreducible component of $\conchsimple$.
\begin{itemize}
\item[(1)] If   $\conchsimple$ is reducible,  the
restricted map $\pi_{2}\vert _{\pi_{1}^{-1}({\cal M})}:
\pi_{1}^{-1}({\cal M})
 \longrightarrow \cc$ is birational.
\item[(2)] The restricted map
$\pi_{1}\vert_{\pi_{1}^{-1}({\cal M})}: \pi_{1}^{-1}({\cal
M})\longrightarrow {\cal M}$ is birational iff ${\cal M}$ is simple.
\end{itemize}
\end{lemma}

\section{Rational Conchoids}\label{sec-rational-conch}
We know that conchoids are either irreducible or with two components
(see Theorem \ref{th-comp}). In this section, we characterize the
conchoids having all their components rational. We see that these
conchoids can only be generated by rational curves. Moreover, we
characterize the cases where the conchoid is rational or it is
reducible with the two components rational. As a consequence, we
prove that the conchoid of the irreducible conics, from a focus on
the conic, are rational;
 in particular all {\it  Lima\c{c}ons of Pascal}  are rational.
We also see that all conchoids of {\it Nicomedes} are rational.

\para

Let $\cP(t)$ be a rational parametrization of $\cc$. Taking into
account the definition of the incidence variety $\bbsimple$, one has
that
\[ \left(\cP(t)+ \frac{d}{\pm \|\cP(t)-A\|}(\cP(t)-A),\cP(t),1+\frac{d}{\pm \|\cP(t)-A\|}\right)\in \bbsimple. \]
So, ${\cal T}^{\pm}(t)=\cP(t)+ \frac{d}{\pm
\|\cP(t)-A\|}(\cP(t)-A)\in \conchsimple$. Therefore, if $\pm
\|\cP(t)-A\|\in \K(t)$,  ${\cal T}^{\pm}(t)$ parametrizes the
components of $\conchsimple$. This motivates the next definition.


\begin{definition}\label{def-rpn} We say that a  parametrization
 ${\cal P}(t) \in
\K(t)^{2}$  {\sf is at rational distance to the focus} if
$\|\cP(t)-A\|^2= m(t)^2,$with $m(t) \in \K(t)$. For short, we
express this fact  saying that ${\cal P}(t)$ is \rdf or $A$-\rdf if
we need to specify the focus.
\end{definition}


\begin{remark}\label{Remark-a-def-rpn}
Note that:
\begin{enumerate}
\item The notion of \rdf depends on the focus. For instance, if  $\cP(t)=(P_{1}, P_{2})=(t,t^2)$ then  $$P_{1}^{2}+P_{2}^{2}=t^2(t^2+1), \hspace{2 mm}  \mbox{and  $P_{1}^{2}+(P_{2}-\frac{1}{4})^{2}=
\frac{1}{16}(1+4t^2)^2$}.$$ So $\cP(t)$ is not  $(0,0)$-\rdf but is
$(0,1/4)$-\rdf.
\item
 If $\cP(t)$ is $A$-\rdf, every re-parametrization of $\cP(t)$ is also $A$-\rdf. However, if can happen
that a re-parametrization of a non \rdf parametrization is \rdf. For
instance, as we have seen above $(t,t^2)$ is not $(0,0)$-\rdf but
$(\frac{2t}{t^2-1},\frac{4t^2}{(t^2-1)^2})$ is $(0,0)$-\rdf since
 \[ \left(\frac{2t}{t^2-1}\right)^2+\left(\frac{2t}{t^2-1}\right)^4=\frac{4 t^2 (t^2+1)^2}{(t^2-1)^4}. \]
 So, we have that if $\cc$ has a proper \rdf parametrization then all the parametrizations of $\cc$ are \rdf with
 respect to the same focus. Nevertheless, it might happen that $\cc$ does not have  proper \rdf
 parametrizations but has non-proper  \rdf pa\-ra\-me\-tri\-zations.
\end{enumerate}
   \end{remark}


Checking whether a given parametrization is \rdf is easy. However
deciding, and actually computing, the existence of \rdf
reparametrizations of non \rdf parametrizations is not so direct.
For dealing with this, we introduce the next notion.


\begin{definition}\label{def-gp} Let $\cP(t) \in
\K(t)^{2}$ be a  rational
 parametrization of $\cc$.
 We define the {\sf reparametrizing
curve of  ${\cal P}(t)$}, and we denote it by $\gp$, as the curve
generated by the primitive part with respect to $x_{2}$ of the
numerator of $\mathfrak{b}((-2x_{2},x_{2}^{2}-1),\cP(x_1)-A)$.
\end{definition}


\begin{remark}\label{Remark-a-def-gp}
We observe that:
\begin{enumerate}
\item $\gp$ does not depend on the representatives of the rational functions in $\cP(t)$.
 \item
The defining polynomial  of $\gp$ has degree 2 w.r.t. $x_2$ and it
is primitive w.r.t. $x_2$. So,  if $\gp$ is reducible then it has
two factors, both depending linearly on $x_2$.
\item Let $\cP(t), {\cal Q}(t)$ be  parametrizations of $\cc$, and $\varphi(t)\in\K(t)$ such that  ${\cal Q}(t)=\cP(\varphi(t))$. Let $$M_1=\mathfrak{b}((-2x_{2},x_{2}^{2}-1),{\cal P}(x_1)-A),\,\,M_2=\mathfrak{b}((-2x_{2},x_{2}^{2}-1),{\cal Q}(x_1)-A).$$
    Then, $  M_1(\varphi(x_1),x_2)=M_2(x_1,x_2).$
\qed
\end{enumerate}   \end{remark}


The following theorem characterizes the conchoids, having  all the
components rational, by means of the notions of \rdf and
reparametrizing curve. In fact, we show that conchoids having all
their components rational can only be generated by rational curves;
indeed  iff the base curve is rational and has  \rdf
parametrizations.


\begin{theorem}\label{th-racional}
  The following statements are equivalent:
\begin{itemize}
\item[(1)]  ${\cal C}$ is rational and has an \rdf parametrization.
\item[(2)] $\conchsimple$ has at least one rational
simple component.
\item[(3)]  There exists a  proper parametrization  of ${\cal C}$ which reparametrizing curve  has at least one rational
 component.
 \item[(4)] The reparametrizing curve of  every proper parametrization  of ${\cal C}$ has at least one rational
 component.
\item[(5)]  All the components of $\conchsimple$
 are rational.
\end{itemize}
\end{theorem}

\noindent  {\bf Proof.} We prove that all the statements are
equivalent to (1). To prove that (2) implies (1), let ${\cal M}$ be
a rational simple component of $\conchsimple$
 parametrized  by ${\cal R}(t)=(R_1(t),R_2(t))$. We consider the
diagram:
\begin{center}
\[ \begin{array}{c}
\Gamma=\pi_{1}^{-1}({\cal M}) \subset \mathfrak{B}(\cc)
  \subset {\Bbb  K}^{2} \times {\Bbb  K}^{2} \times {\Bbb  K} \\
  \mbox{\hspace*{3.8 cm} $\begin{array}{c}
\vspace*{-2 mm}
\tilde{\pi}_{1}=\pi_{1}\vert_{\pi_{1}^{-1}({\cal M})} \\
\\
\end{array}$ \hspace{-2 mm} {\Huge $\swarrow$}
} \hspace{2 mm} \mbox{ {\Huge $\searrow$} $\begin{array}{c}
\vspace*{-2 mm}
\tilde{\pi}_{2}=\pi_{2}\vert_{\pi_{1}^{-1}({\cal M})} \\
\\
\end{array}$}   \\
\hspace*{3 cm} {\cal M} \subset \conchsimple \subset
\K^{2} \hspace{1 cm}    {\cal C} \subset \K^{2} \\
 \mbox{ \hspace*{-1.3cm} {\Huge $\uparrow$} ${\cal R}$} \\
\mbox{\hspace*{-1.6 cm}  $\K$}
\end{array}  \]
\end{center}
Since ${\cal M}$ is simple,  by Lemma \ref{lemma-pi-birracional},
$\tilde{\pi}_{1}$ is birational. So, ${\cal
Q}(t)=\tilde{\pi}_{2}(\tilde{\pi}_{1}^{-1}({\cal R}(t)))$
parametrizes ${\cal C}$. Let us see that ${\cal
Q}(t):=(Q_1(t),Q_2(t))$ is \rdf. By construction,
$\tilde{\pi}_{1}^{-1}({\cal R}(t))=({\cal R}(t), {\cal
Q}(t),\lambda) \in \bbsimple$, where
$\lambda=(R_1-a)/(Q_1-a)=(R_2-b)/(Q_2-b)$ and $\|{\cal R}(t)-{\cal
Q}(t)\|^2=d^2$. Note that $\cc$ is not a line passing through $A$,
and hence $Q_1\neq a$, $Q_2\neq b$. Moreover, $\lambda \neq 1$,
since otherwise  ${\cal R}(t)={\cal Q}( t)$ that yields to $d=0$. So
${\cal Q}(t)-A= ({\cal Q}(t)-{\cal R}(t))+({\cal R}(t)-A)=({\cal
Q}(t)-{\cal R}(t))+\lambda({\cal Q}(t)-A)$, and hence   $\|{\cal
Q}(t)-A\|^2=d^2/(\lambda-1)^2.$

  In order to prove that (1) implies
(2), let  ${\cal P}(t)$ be an \rdf parametrization of $\cc$. Let
$\|\cP(t)-A\|^2=m(t)^2$.  Then
\[ \left(\cP(t)\pm \frac{d}{m(t)}(\cP(t)-A),\cP(t),1\pm\frac{d}{m(t)}\right)\in \bbsimple. \]
Moreover, since $\cP(t)$ generates a dense subset of $\cc$, by Lemma
3 in \cite{SS08}, $\cP(t)\pm \frac{d}{m(t)}(\cP(t)-A)$ generates a
dense in $\conchsimple$. So, all components of $\conchsimple$ are
rational. Now, the result follows from  Theorem
\ref{th-properties-simple-special}.

To see that (1) implies (3), let ${\cal P}(t)=(P_{1}(t),P_{2}(t))$
be an \rdf  pa\-ra\-me\-tri\-za\-tion
 of $\cc$, and $\| \cP(t)-A\|^2=m(t)^2$. Then, $1/m(t) (\cP(t)-A)$ parametrizes the circle
  $x_{1}^{2}+x_{2}^{2}=1$. Since
${\cal R}(t)=(\frac{t^{2}-1} {t^{2}+1}, \frac{2\,t}{t^2+1})$ is a
proper parametrization of the circle, it holds that there exists
$\phi(t) \in \K(t)$ such that ${\cal R}(\phi(t))=1/m(t)(\cP(t)-A)$.
This implies that
$\mathfrak{b}((-2\phi(t),\phi(t)^2-1),\cP(t)-A)=0.$ Therefore $(t,
\phi(t))$ parametrizes one component of $\gp$.

To prove that (3) implies (1), let $(\phi_{1}(t),\phi_{2}(t))$ be a
parametrization of one component of $\gp$, where
$\cP(t)=(P_1(t),P_2(t))$ is a proper parametrization of $\cc$. Then,
$\mathfrak{b}((-2\phi_2(t),\phi_{2}(t)^2-1),\cP(\phi_1(t))-A)=0$.
Note that $\phi_{2}$ is not identically zero since otherwise it
would imply that $P_{2}(\phi_{1})=b$ and $\cc$ is not a line passing
through the focus. Then, it follows that $\cP(\phi_1(t))$ is \rdf;
indeed
\[ \| \cP(\phi_1(t))-A\|^2=\frac{(\phi_{2}(t)^{2}+1)^2}{(2\phi_2(t))^2}(P_2(\phi_1(t))-b)^2. \]

 Trivially (4)
implies (3). In order to prove that (3) implies (4),  let ${\cal
P}(t)$ and ${\cal Q}(t)$ be two proper parametrizations of ${\cal
C}$, such that $\gp$ has at least one rational component. Let
$(\phi_1(t),\phi_2(t))$ be a parametrization of one component of
$\gp$. Note that, because of Remark \ref{Remark-a-def-gp} (2),
$\phi_1(t)$ is not constant. Since both parametrizations are proper,
there exists an invertible $\varphi \in \K(t)$ such that ${\cal
Q}(t)={\cal P}(\varphi(t))$.  Let
$M_1(x_1,x_2)=\mathfrak{b}((-2x_2,x_{2}^2-1),\cP(x_1)-A)$ and
$M_2(x_1,x_2)=\mathfrak{b}((-2x_2,x_{2}^2-1),{\cal Q}(x_1)-A)$.Let
$D_i$ be the denominator of $M_i$ and  let $C_i,H_i$ be,
respectively, the content and primitive part w.r.t. $x_2$ of the
numerator of $M_i$.  Then, by Remark \ref{Remark-a-def-gp} (3),
\[ C_1(x_1)H_1(x_1,x_2)D_2(\varphi^{-1}(x_1))=D_1(x_1) C_2(\varphi^{-1}(x_1))H_2(\varphi^{-1}(x_1),x_2). \]
So, $D_1(\phi_1)
C_2(\varphi^{-1}(\phi_1))H_2(\varphi^{-1}(\phi_1),\phi_2)=0$. Since
$\phi_1\not\in \K$, then $\varphi^{-1}(\phi_1(t))\not\in \K$. Since
 $D_1,C_2$ are non-zero univariate polynomials,  $D_1(\phi_1) C_2(\varphi^{-1}(\phi_1))\neq 0$. Therefore, $H_2(\varphi^{-1}(\phi_1),\phi_2)=0$. Hence, $(\varphi^{-1}(\phi_1),\phi_2)$ parametrizes a component of $\mathfrak{G}({\cal Q})$.
Therefore one concludes (4). For the implication of (1) implies (5)
see the proof of (1) implies (2). Furthermore, if (5) holds, then
$\conchsimple$ has at least one rational simple component, and by
(2) one concludes (1). \qed


\begin{remark}\label{remark-a-th-racional}
Theorem \ref{th-racional} implies that:
\begin{enumerate}
\item Conchoids with all their
components rational can only be generated by   rational curves.
\item The rationality of all the components of the conchoid does depend on the base curve and the focus, but not on the distance.
\end{enumerate}
\end{remark}


\begin{corollary}\label{cor-types-of-conch}
The conchoid of a rational curve is either rational, or  it is
reducible with two rational components, or it is irreducible but
non-rational.
\end{corollary}

\noindent  {\bf Proof.} Let $\cc$ be the base curve, and let
$\conchsimple$ be reducible. Then, by Corollary 3 in \cite{SS08} and
Theorem \ref{th-properties-simple-special}, at least one  conchoid
simple component is rational. Now the corollary follows from
Theorem \ref{th-racional}.  \qed


\begin{remark}\label{remark-a-th-racional}
If $\cc$ is non-rational, it might happen that its conchoid is
reducible with two non-rational components or with one component
non-rational and the other rational. For instance, if $\cc$ is the
curve defined by the polynomial

\vspace{2 mm}
 \noindent
$f(y_1,y_2)=81+162y_2^{2}y_1^6+1521y_1^4y_2^2+972y_1y_2^4+162y_1^2y_2^4-1944y_1^3y_2^3-1944y_1^5y2+
864y_1^3y_2^2-2898y_1^2y_2^3-8730y_1^4y_2-1404y_1^2-1458y_2+108y_1-324y_1^4y_2^3+81y_2^4y_1^4-
17388y_1^3y_2+972y_1^5y_2^2-162y_2^5y_1^2-162y_1^6y_2-13122y_1^2y_2+6480y_2^2-162y_2^3-972y_2y_1+8694y_2^2y_1
+4203y_1^2y_2^2+1449y_2^4-1932y_1^3+5590y_1^4+972y_1^7+81y_1^8+8532y_1^5+4356y_1^6+81y_2^6
$

\vspace{2 mm}
 \noindent then the conchoid of  $\cc$  from the focus $A=(-3,0)$ and
distance $d=1/3$, has two components  defined by the factors

\vspace{2 mm}
 \noindent
$(-x_2+x_1^2)\,(1296+1728x_1-5832x_2+6237x_2^2-3888x_2*x_1+3690x_1^2x_2^2-4968x_1^2-13122x_1^2
x_2-7584x_1^3+2689x_1^4-162x_1^6x_2+1422x_2^4+972x_1^7-162x_2^5x_1^2+972x_1^5x_2^2-17064x_1^3x_2-
648x_2^3+81x_2^4x_1^4-8676x_1^4x_2+8532x_1x_2^2+81x_1^8-2844x_1^2x_2^3+7884x1^5+540x_1^3x_2^2-
1944x_1^5x_2+4302x_1^6-1944x_1^3x_2^3+162x_1^2x_2^4+972x_1x_2^4+1467x_1^4x_2^2+162x_2^2x_1^6+81x_2^6
-324x_1^4x_2^3) $,

\vspace{2 mm}
 \noindent one of them is the parabola (rational)
and the other is a non-rational curve of genus $1$.

    Note that the reason is that simple components of reducible conchoids are birationally equivalent to $\cc$ (see Corollary  3 in \cite{SS08}). Also, take into account that for almost all values of $d$ the conchoid has all the components simple (see Theorem 4 in \cite{SS08}).
\end{remark}


\begin{corollary}\label{corolario-a-th-racional}
Let   ${\cal P}(t)$ be a parametrization of ${\cal C}$ such that
$\gp$ has at least one rational component ${\cal M}$, and let
$(\phi_{1}(t),
 \phi_{2}(t))$ be a parametrization of ${\cal M}$. Then   ${\cal P} (\phi_{1}(t))$ is \rdf.
\end{corollary}

\noindent  {\bf Proof.}
 It follows from the proof of (3) implies (1) in Theorem \ref{th-racional}. \qed


\begin{corollary}
\label{P41} Let $\cP(t)$ be a proper parametrization of ${\cal C}$.
Then, the following statements are equivalent:
\begin{itemize}
\item[(1)]  All the components of $\conchsimple$ are rational.
\item[(2)] There exists $\varphi \in \K(t)$
 of degree at most two such that ${\cal
P}(\varphi(t))$ is \rdf.
\end{itemize}
\end{corollary}

\noindent  {\bf Proof.} (1) implies (2) follows from Theorem
\ref{th-racional}.  (2) implies (1) follows from Theorem
\ref{th-racional}, from Corollary \ref{corolario-a-th-racional}, and
using that the partial degree of $\gp$ w.r.t. $x_2$ is 2 (see
Theorem 4.21 in \cite{libro}). \qed

\para

In the sequel, we analyze the case of conchoids of rational curves,
characterizing the rational conchoid and conchoids with two rational
components; we refer to this case as {\sf double rationality}.


\begin{lemma}
\label{lemma-gp} $\gp$ is reducible if and only if ${\cal P}(t)$ is
\rdf.
\end{lemma}

\noindent  {\bf Proof.} Let $H$ be the primitive part w.r.t. $x_2$,
and $M(x_1)$ the content w.r.t. $x_2$ of the numerator of
$\mathfrak{b}((-2x_{2},x_{2}^{2}-1),\cP(x_1)-A)$. All factors of $H$
depend
 on $x_{2}$. Thus,
 $\gp$ is reducible
 if and only if $H$ has  two factors depending on $x_{2}$, or
equivalently, the discriminant $\Delta_{H}$ w.r.t. $x_{2}$ is the
square of a polynomial. Therefore, since
$M(x_1)^2\Delta_{H}=4\,\|{\cal P}(t)-A \|^2$,
 one has that $\gp$
is reducible if and only if $\cP(t)$ is \rdf. \qed




\begin{theorem}
\label{theorem-reducible} {\sf (Characterization of double rational
conchoids)} Let ${\cal C}$ be rational. The following statement are
equivalent:
\begin{itemize}
\item[(1)] $\conchsimple$ is reducible.
\item[(2)]  $\conchsimple$ has exactly two components and they are
rational.
\item[(3)] There exists an \rdf  proper   parametrization
 of $\cc$.
\item[(4)]  Every proper parametrization  of $\cc$ is \rdf.
\item[(5)]  There exists a proper parametrization of $\cc$ which reparametrizing curve is reducible.
\item[(6)]  The reparametrizing curve of  every proper parametrization of $\cc$  is reducible.
\end{itemize}
\end{theorem}

\noindent  {\bf Proof.}
By Corollary \ref{cor-types-of-conch}, (1) implies (2).  (2) implies
(1) trivially. In order to prove that (2) implies (3),
 let ${\cal R}(t)$ be a proper  parametrization of
a simple component  ${\cal M}$ of $\conchsimple$. We consider the
diagram used in the proof of Theorem \ref{th-racional}. By Lemma
\ref{lemma-pi-birracional},  $\tilde{\pi}_{2} \circ
\tilde{\pi}_{1}^{-1} \circ {\cal R}: \K \longrightarrow \cc$ is
birational.
 Therefore, ${\cal Q}(t)=\pi_{2}(\pi_{1}^{-1}({\cal R}(t)))$ is
 a proper parametrization of ${\cal C}$. Furthermore,
reasoning as in the proof of  ``(2) implies (1)", in Theorem
\ref{th-racional}, one has that ${\cal Q}(t)$ is \rdf. (3) implies
(4) follows from the Remark \ref{Remark-a-def-rpn}.

In order to  see that (4) implies (2), let ${\cal
P}(t)=(P_1(t),P_2(t))$ be an \rdf proper
 parametrization of ${\cal C}$. Reasoning as in the proof of Theorem \ref{th-racional}, we get that
  $${\cal R}^{\pm}(t)=(R_{1}^{\pm}(t),R_{2}^{\pm}(t)):={\cal P}(t) \pm d \, \frac{{\cal P}(t)-A}{\| {\cal
P}(t)-A\|}$$ parametrizes  all components of $\conchsimple$. So, it
only remains to prove that $\conchsimple$ is reducible. Let us
assume that it is  irreducible. Then, by Theorem
\ref{th-properties-simple-special}, $\conchsimple$ is simple, and,
by Lemma \ref{lemma-pi-birracional} (2), $\bbsimple$ is irreducible.
Moreover,
\[ {\cal M}^{\pm}(t)=\left({\cal R}^{\pm}(t),{\cal P}(t),\frac{R_{1}^{\pm}(t)-a}{P_1(t)-a}\right) \]
are two rational parametrizations of $\bbsimple$. Moreover, since
$\cP(t)$ is proper then $\K(t)=\K(\cP(t))\subset \K({\cal
M}^{\pm}(t))\subset \K(t).$ So, each ${\cal M}^{\pm}(t)$ is proper.
So, there exists a linear rational function $\varphi(t)$ such that
${\cal M}^{+}(\varphi(t))={\cal M}^{-}(t)$. Thus, $\varphi(t)=t$
and, since $d\neq 0$,  $\cP(t)=A$ which is a contradiction.

Applying Lemma \ref{lemma-gp} one has  that (4) implies (5). The
implication ``(5) implies (6)" follows from  Lemma \ref{lemma-gp}
and  Remark
 \ref{Remark-a-def-rpn}.
Finally, ``(6) implies (4)" follows directly from Lemma
\ref{lemma-gp}. \qed


\begin{theorem}
\label{theorem-irreducible} {\sf (Characterization of  rational
conchoids)} Let ${\cal C}$ be rational. The following statement are
equivalent:
\begin{itemize}
\item[(1)]  $\conchsimple$ is  rational.
\item[(2)] There exists a proper parametrization
of ${\cal C}$ which reparametrizing curve is rational.
\item[(3)]  The reparametrizing curve of  every proper parametrization
of ${\cal C}$ is rational.
\end{itemize}
\end{theorem}


\noindent  {\bf Proof.}
 Let  $\conchsimple$ be
rational. By Theorem \ref{th-racional}, there exists a proper
pa\-ra\-me\-tri\-za\-tion ${\cal P}(t)$ of ${\cal C}$ such that
 $\gp$ has at least one rational component; say
${\cal M}$. Furthermore, by Theorem \ref{theorem-reducible}, ${\cal
P}(t)$ is not \rdf. Thus, by Lemma \ref{lemma-gp}, $\gp$ is
irreducible, and hence rational. So, (1) implies (2).

We prove that (2) implies (3). Let $\cP(t)$ proper such that $\gp$
is rational and let ${\cal Q}(t)$ be another proper parametrization
of $\cc$. Since $\gp$ is irreducible, by Lemma \ref{lemma-gp},
$\cP(t)$ is not \rdf. By Remark \ref{Remark-a-def-rpn} ${\cal Q}(t)$
is not \rdf. So,  by Lemma \ref{lemma-gp}, $\mathfrak{G}({\cal Q})$
is irreducible. Now, the result follows from Theorem
\ref{th-racional}.

Finally, we prove that (3) implies (1).  Let $\gp$ be rational, with
$\cP(t)$ proper. By Lemma \ref{lemma-gp}, $\cP(t)$ is not \rdf.
Thus, by Theorem \ref{theorem-reducible}, $\conchsimple$ is
irreducible. Now, the result follows from Theorem \ref{th-racional}.
\qed

\para

We apply these results to the case of conchoids of conics with the
focus on the conic (in particular to {\it  Lima\c{c}ons of Pascal}),
and to the case of conchoids of {\it Nicomedes}.

\begin{lemma}\label{lemma-partial-degree}
Let $\cP(t)=(p_{1}(t)/p(t),p_2(t)/p(t))$ be a proper parametrization
of  $\cc$ with $\gcd(p_1,p_2,p)=1$ and $\deg_t(p_i/p)\leq 2$ for
$i=1,2$. If $A\in \cc$, then  $\conchsimple$ is rational.
\end{lemma}

\noindent {\bf Proof.} The defining polynomial $g(x_1,x_2)$ of $\gp$
is the primitive part w.r.t. $x_2$ of $K(x_1,x_2)=-2x_2 (p_1(x_1)-a
p(x_1))+(x_{2}^{2}-1)(p_2(x_1)-b p(x_1)).$  Moreover, the content
$C(x_1)$ of $K$ w.r.t. $x_2$ is $\gcd(p_1(x_1)-ap(x_1),p_2(x_1)-b
p(x_1))$. First, we observe that $\deg_{x_1}(g)>0$. Indeed, if it is
zero, it implies that there exist $\lambda, \mu \in \K$ such that
$\cP(t)=(a+\lambda C(t)/p(t),b+\mu C(t)/p(t))$ and, hence, $\cc$
would be a line passing through the focus.

Let us assume that $A=(a,b)$ is reachable by $\cP(t)$; say
$\cP(t_0)=A$. Then $x_1-t_0$ divides $C(x_1)$, and hence
$\deg_{x_1}(g)=1$. So  $\gp$ is rational and, by Theorem
\ref{theorem-irreducible}, $\conchsimple$ is rational.  Now, if $A$
is not reachable by $\cP(t)$. Then, for $i=1,2$, $\deg(p_i)\leq
\deg(p)$ (see Section 6.3. in \cite{libro}). Say
$$p_i(x_1)=a_{i,n}x_{1}^{n}+\cdots +a_{i,0},\,\, p(x_1)=b_n x_{1}^{n}+\cdots+b_0,$$
where $a_{i,n}$ might vanishes. Then, $A=(a_{1,n}/b_n,a_{2,n}/b_n)$
(see Section 6.3. in \cite{libro}). So, $\deg_{x_1}(g)=1$. Thus,
reasoning as above we get the result. \qed

%

\para

Now, taking into account the  parametrizations of the irreducible
conics, by Lemma \ref{lemma-partial-degree}, one deduces the
following result (see also Examples \ref{ex-parabola},
\ref{ex-ellipse}, \ref{ex-pascal}, \ref{ex-hyperbola}).


\begin{corollary}\label{corollary-conics}
Let $\cc$ be an irreducible conic, and let $A\in \cc$, then
$\conchsimple$ is rational.
\end{corollary}

\begin{corollary}\label{corollary-pascal}
{\it  Lima\c{c}ons of Pascal} are rational.
\end{corollary}


\begin{remark}
In general it is not true that if the focus is on the curve,  the
conchoid is rational. For instance, let $\cc$ be the curve defined
by $y_{1}^{3}y_{2}=1$, $\cP(t)=(1/t,t^3)$, and $A=(1,1)\in \cc$.
Then, $\gp$ is defined by
$x_{2}^2x_{1}^3-x_{1}^3+x_{1}^2x_{2}^2-x_{1}^2+x_{1}x_{2}^2-x_{1}+2x_{2}$
and its genus is 2.
\end{remark}

Finally we analyze the conchoids of {\it Nicomedes} (see also
Example \ref{ex-nicomedes}).

\begin{corollary}\label{corollary-nicomedes}
Conchoids of {\it Nicomedes} are rational.
\end{corollary}

\noindent {\bf Proof.} Conchoids of {\it Nicomedes} appear when
$\cc$ is a line and $A\not\in \cc$. Let
$\cP(t)=(p_1(t),p_2(t))=(a_1+\lambda_1 t, a_2+\lambda_2 t)$.  The
defining polynomial  of $\gp$ is $g(x_1,x_2)=-2x_2
(p_1(x_1)-a)+(x_{2}^{2}-1)(p_2(x_1)-b).$ Note that, since
$A\not\in\cc$, $g$ is primitive w.r.t. $x_2$. Now the result follows
from Theorem \ref{theorem-irreducible} and noting that
$\deg_{x_1}(g)=1$.  \qed

\section{Parametrization of Conchoids.}\label{sec-param-conch}

In this section we apply the results in Section
\ref{sec-rational-conch} to derive an algorithm to check the
rationality of the components of a conchoid and, in the affirmative
case, to parametrize them. For this purpose, in the sequel, let
$\cc$ be rational and $\cP(t)$ be a proper  parametrization of
$\cc$. We also assume that the focus $A$ is fixed. However, we
consider $d$ generic. Recall that we have assume that $\cc$ is not a
line through the focus neither a circle centered at the focus and
radius $d$; nevertheless,  observe that the problem for these two
excluded cases is trivial.

\para

First, we check whether $\cP(t)$ is \rdf; equivalently one can check
whether $\gp$ is reducible. If so,  by Theorem
\ref{theorem-reducible}, $\conchsimple$ is double rational and
$$\cP(t)+ \frac{d}{\pm \|\cP(t)-A\|}(\cP(t)-A)$$
parametrizes the two components. If $\cP(t)$ is not \rdf, we check
whether $\gp$ is rational. If it is not rational, by Theorem
\ref{theorem-irreducible}, $\conchsimple$ is not rational. If $\gp$
is rational, by Theorem \ref{theorem-irreducible},  $\conchsimple$
is rational. In order to parametrize $\conchsimple$, we get a proper
parametrization $(\phi_1(t),\phi_2(t))$ of $\gp$ (see \cite{libro}
for this). Then, by Corollary \ref{corolario-a-th-racional}, ${\cal
Q}(t)=\cP(\phi_1(t))$ is \rdf. Therefore, any of the
parametrizations
$${\cal Q}(t)+ \frac{d}{\pm \|{\cal Q}(t)-A\|}({\cal Q}(t)-A)$$
paramatrizes $\conchsimple$. Summarizing we get the following
procedure:
\begin{enumerate}
\item Compute the primitive part $g(x_1,x_2)$ w.r.t. $x_2$ of the numerator of \\ $\mathfrak{b}((-2x_2,x_{2}^{2}-1),\cP(x_1)-A)$.
\item If $g$ is reducible {\sf return} that $\conchsimple$ is double rational and that \\  $\cP(t)+ \frac{d}{\pm \|\cP(t)-A\|}(\cP(t)-A)$ parametrize the two components.
\item Check whether the genus of $\gp$ is zero. If not {\sf return} that $\conchsimple$ is not rational.
\item Compute a proper parametrization $(\phi_1(t),\phi_2(t))$ of $\gp$ and {\sf return} that $\conchsimple$ is  rational and that  $\cP(\phi_1(t))+ \frac{d}{\pm \|\cP(\phi_1(t))-A\|}(\cP(\phi_1(t))-A)$ pa\-ra\-me\-trizes $\conchsimple$.
\end{enumerate}

We illustrate the algorithm by means of some examples.


\begin{example}\label{ex-parabola} {\sf (Conchoid of Parabolas)}
Let $\cc$ be the parabola defined by $f(y_1,y_2)=y_2-\mu_1
y_{1}^{2}+\mu_2 y_1+\mu_3$, with $\mu_1\neq 0$. We consider the
proper parametrization $\cP(t)=(t,\mu_1 t^2+\mu_2 t+\mu_3),$ and the
focus $A=(\lambda, \mu_1 \lambda^2+\mu_2 \lambda +\mu_3)$ being any
point on $\cc$. By Corollary \ref{corollary-conics}, we know that
$\conchsimple$ is rational. Here, we indeed compute a
parametrization. The polynomial $g$ defining $\gp$ is irreducible:
$$ g(x_1,x_2)=\mu_1x_1x_{2}^{2}-\mu_1 x_{1}+ \mu_2 x_{2}^{2}-\mu_2-2x_{2}+\lambda \mu_1 x_{2}^{2}-\lambda \mu_1.$$
Moreover $\gp$ is rational and can be parametrized as (recall that
$\mu_1\neq 0$)
\[ \phi(t)=(\phi_1(t),\phi_2(t))=\left(-{\frac {{t}^{2}\mu_2-\mu_2-2\,t+\lambda\,{t}^{2}\mu_1-\lambda\,\mu_1}{\mu_1 \left( {t}^{
2}-1 \right) }},t\right). \] Therefore, $ {\cal
Q}(t)=(\phi_1(t),\mu_1 \phi_1(t)^2+\mu_2\phi_1(t)+\mu_3) $ is \rdf
and ${\cal Q}(t)+ \frac{d}{\pm \|{\cal Q}(t)-A\|}({\cal Q}(t)-A)$
parametrizes $\conchsimple$.
\end{example}


\begin{example}\label{ex-ellipse} {\sf (Conchoid of Ellipses)}
Let $\cc$ be the ellipse  defined by
$$f(y_1,y_2)= \frac{y_{1}^{2}}{r_{1}^{2}}+\frac{y_{2}^{2}}{r_{2}^{2}}-1,$$ with $r_1 r_2 \neq 0$. We consider the proper parametrization
$$\cP(t)=\left({\frac {2\,r_{{1}}t}{{t}^{2}+1}},{\frac {r_{{2}} \left( {t}^{2}-1
 \right) }{{t}^{2}+1}}\right),$$ and the focus $A=\cP(\lambda)$ being a point on $\cc$. By Corollary \ref{corollary-conics}, we know that $\conchsimple$ is rational. Here, we indeed compute a parametrization. The polynomial $g$,  defining $\gp$, is irreducible:
$$ g(x_1,x_2)=2\,x_{{1}}x_{{2}}r_{{1}}\lambda+x_{{1}}r_{{2}}{x_{{2}}}^{2}-r_{{2}}x_{
{1}}-2\,x_{{2}}r_{{1}}+\lambda\,r_{{2}}{x_{{2}}}^{2}-r_{{2}}\lambda.$$
Moreover $\gp$ is rational and can be parametrized as
\[ \phi(t)=(\phi_1(t),\phi_2(t))=\left(-{\frac {-2\,r_{{1}}t+\lambda\,r_{{2}}{t}^{2}-r_{{2}}\lambda}{2\,t r_{
{1}}\lambda+r_{{2}}{t}^{2}-r_{{2}}}},t\right). \] Therefore, ${\cal
Q}(t)=\cP(\phi_1(t))$ is \rdf and ${\cal Q}(t)+ \frac{d}{\pm \|{\cal
Q}(t)-A\|}({\cal Q}(t)-A)$ parametrizes $\conchsimple$.
\end{example}

\begin{example}\label{ex-pascal} {\sf (Lima\c{c}on of Pascal)}
Taking $r_1=r_2\neq 0$ in Example \ref{ex-ellipse}, we get a
parametrization of the {\it  Lima\c{c}ons of Pascal}
\end{example}


\begin{example}\label{ex-hyperbola} {\sf (Conchoid of Hyperbolas)}
Let $\cc$ be the hyperbola  defined by
$$f(y_1,y_2)= \frac{y_{1}^{2}}{r_{1}^{2}}-\frac {y_{2}^{2}}{r_{2}^{2}}-1,
$$ with $r_1 r_2 \neq 0$. We consider the proper parametrization
$$\cP(t)=\left(-\frac{r_{1} \left(r_{1}^{2}+r_{2}^{2}{t}^{2} \right) }{
-r_{1}^{2}+r_{2}^{2}{t}^{2}},\frac{2r_{2}^{2}r_{1}
t}{-r_{1}^{2}+r_{2}^{2}t^{2}}\right),$$ and the focus
$A=\cP(\lambda)$ being a point on $\cc$. By Corollary
\ref{corollary-conics}, we know that $\conchsimple$ is rational.
Here, we indeed compute a parametrization. The polynomial $g$,
defining $\gp$, is irreducible:
$$ g(x_1,x_2)=-x_{1}x_{2}^{2}\lambda\,r_{2}^{2}-2\,x_{1}x_{2}r_{1}^{2}+\lambda
\,r_{2}^{2}x_{1}-x_{2}^{2}r_{1}^{2}+r_{1}^{2}-2\,\lambda\,x_{2}r_{1}^{2}.$$
Moreover $\gp$ is rational and can be parametrized as
\[ \phi(t)=(\phi_1(t),\phi_2(t))=\left(  -\frac{r_{1}^{2} \left( {t}^{2}-1+2\,\lambda\,t \right) }{
\lambda\,r_{2}^{2}{t}^{2}+2\,r_{1}^{2}t-\lambda\,r_{2}^{2}
},t\right). \] Therefore, ${\cal Q}(t)=\cP(\phi_1(t))$ is \rdf and
${\cal Q}(t)+ \frac{d}{\pm \|{\cal Q}(t)-A\|}({\cal Q}(t)-A)$
parametrizes $\conchsimple$.
\end{example}

\begin{example}\label{ex-nicomedes} {\sf (Conchoid of Nicomedes)}
Let $\cc$  be the line parametrized by
$$\cP(t)=(a_1+t \lambda_1, a_2+t \lambda_2),$$ and the focus $A=(a,b)\not \in \cc$. Then, $\conchsimple$ is the conchoid of {\it Nicomedes}. By Corollary \ref{corollary-conics}, we know that $\conchsimple$ is rational. Here, we indeed compute a parametrization. The polynomial $g$,  defining $\gp$, is irreducible because $A\not\in \cc$:
$$ g(x_1,x_2)=-2 x_2 a_1-2 x_2 x_1 \lambda_1+2 x_2 a+x_{2}^2a_2+x_{2}^{2}x_{1}\lambda_{2}-x_{2}^{2} b-a_2-x_1\lambda_{2}+b.$$
Moreover $\gp$ is rational and can be parametrized as
\[ \phi(t)=(\phi_1(t),\phi_2(t))=\left( \frac{2 t a_1+t^2 b-2 t a-t^2 a_2-b+a_2}{-2 t\lambda_1-\lambda_2+t^2 \lambda_2},t\right). \]
Therefore, ${\cal Q}(t)=\cP(\phi_1(t))$ is \rdf and ${\cal Q}(t)+
\frac{d}{\pm \|{\cal Q}(t)-A\|}({\cal Q}(t)-A)$ parametrizes
$\conchsimple$.
\end{example}

\section{Detecting focuses to Parametrize Conchoids.}\label{sec-param-focus}

In the previous section we have seen how to decide whether the
conchoid to a rational curve is rational or double rational and, in
the affirmative case, how to parametrize the components of the
conchoid. Nevertheless, in that reasoning the focus is fixed. In
this section, we analyze a slightly different problem. We assume
that we are given a proper parametrization
$$\cP(t)=\left(\frac{p_1(t)}{p(t)},\frac{p_2(t)}{p(t)}\right),$$
where $\gcd(p_1,p_2,p)=1$, of a rational curve $\cc$, and we look
for $A_0\in \K^2$ such that the conchoid $\mathfrak{C}(\cc,A_0,d)$
has all the components rational. We know that this implies that
either $\mathfrak{C}(\cc,A_0,d)$ has two rational components or it
is rational. In the first case we say that the $A_0$ is a {\sf
double rational focus} and, in the second, that $A_0$ is a {\sf
rational focus}. For this purpose, in the sequel, $A=(a,b)$ is
treated generically, and hence $a,b$ are unknowns.

\vspace{2 mm}

\noindent {\sf \underline{Detecting double rational focuses.}} The
strategy is as follows. First we determine a set $\cal F$ in  $\K^2$
containing the possible double rational focuses. Afterwards, we
prove that $\cal F$  is the union of $\cc$ and finitely many lines.
So  all  components of $\cal F$ are rational, and using a
parametrization of each component we determine conditions on the
parameter to get double rational focuses.

\para

More precisely, let $$\Delta_1(a,t)=p_1(t)-a\, p(t),\,\,
\Delta_2(b,t)=p_2(t)-b\, p(t),$$
$$\Sigma_1(a,b,t)=\Delta_1+\sqrt{-1}\, \Delta_2,\,\,
\Sigma_2(a,b,t)=\Delta_1-\sqrt{-1}\, \Delta_2.$$ So,
$$\|\cP(t)-A\|^2=\frac{1}{p(t)^2}(\Delta_{1}^{2}+\Delta_{2}^{2})=\frac{1}{p(t)^2}\Sigma_1
\Sigma_2.$$ Thus, a necessary condition for $A_0=(a_0,b_0)\in\K^2$
to be a double rational focus of $\cP(t)$ is that
$\Sigma_1(a_0,b_0,t) \Sigma_2(a_0,b_0,t)$ is   either constant or it
has multiple roots. Let us see that the first condition cannot
happen; recall that we have excluded the case where $\cc$ is a line
and $A_0\in \cc$.


\begin{lemma}\label{lemma-double-rational-focus}
For every focus  $A_0=(a_0,b_0)\in \K^2$, $\Sigma_1(a_0,b_0,t)
\Sigma_2(a_0,b_0,t)$ is not constant.
 \end{lemma}

\noindent {\bf Proof.} Let $a_0,b_0\in \K$ be such that
$\Sigma_1(a_0,b_0,t) \Sigma_2(a_0,b_0,t)\in \K$. We prove that then
$\cc$ is a line passing through $A_0=(a_0,b_0)$ and we have excluded
this case. It holds that either there exists $i\in \{1,2\}$ such
that $\Sigma_i(a_0,b_0,t)=0$ (say, $i=1$) or both
$\Sigma_i(a_0,b_0,t), i=1,2$, are constant. In the first case,
$p_1(t)/p(t)=a$, and this implies that $\cc$ is the line $y_1=a$. So
$\Sigma_i(a_0,b_0,t)\in \K$ for $i=1,2$.  This implies that
$\Delta_1(a_0,t),\Delta_2(b_0,t)$ are constants. Say
$\Delta_1(a_0,t)=\mu_1,\Delta_2(b_0,t)=\mu_2$, with $\mu_1,\mu_2\in
\K$. However, this implies that $\cc$ is the line
$\mu_2(y_1-a)=\mu_1 (y_1-b)$ that passes through $A_0$. \qed

\para

Therefore, if $A_0$ is a double rational focus then at least one the
following holds: (i) $\Sigma_1(a_0,b_0,t), \Sigma_2(a_0,b_0,t)$ have
a common root or, equivalently, $\Delta_1(a_0,t), \Delta_2(b_0,t)$
have a common root; (ii) $\Sigma_1(a_0,b_0,t)$ has a multiple root;
(iii) $ \Sigma_2(a_0,b_0,t)$ has a multiple root.

\para

Now, let  $R(a,b)$ be the square-free part of the resultant of
$\Delta_{1}(a,t), \Delta_2(b,t)$ w.r.t. $t$, and $D_i(a,b)$ be the
square-free part of the resultant w.r.t. $t$ of $\Sigma_{i}(a,b,t),
\partial \Sigma_{i}(a,b,t)/ \partial t$, respectively; note that
$D_i$ is the square-free part of the  discriminant of $\Sigma_i$
w.r.t $t$ multiplied by the leading coefficient of $\Sigma_i$ w.r.t.
$t$. Then,  the double rational focuses belong to the algebraic set
$\cal F$ in $\K^2$ defined by $R(a,b)D_1(a,b)D_2(a,b)=0$.

\para

By Theorem 4.41 in \cite{libro}, $R$ is  the defining polynomial of
$\cc$. Moreover, since
\[ \Sigma_i(a,b,t)=(p_1(t)\pm \sqrt{-1}\,p_2(t)) - (a \pm \sqrt{-1} \,b) p(t), \]
$D_i(a,b)$ can be expressed as a polynomial in $(a \pm
\sqrt{-1}\,b)$ and hence it is a product of linear factors in $a,b$.
Thus, we have the following proposition.


\begin{proposition}
$\cal F$ decomposes as $\cc$ union finitely many lines.
\end{proposition}


Now,   we take a parametrization ${\cal Q}(h)=(Q_{1}(h),Q_{2}(h))$
of each component of $\cal F$, and we consider the rational function
$\Delta(h,t):=\Delta_{1}(Q_{1}(h),t)^2+\Delta_{2}(Q_{2}(h),t)^{2}$.
Repeating a similar argument as above we determine necessary
conditions on $h$ such that ${\cal Q}(h)$ is a double rational
focus, and a final checking detect the double rational focuses, when
they exist.


\begin{example}\label{example-parabola-double-rat-focus} {\sf (Double rational focuses for the parabola)}
Let $\cc$ be the parabola over $\Bbb C$ parametrized by
$\cP(t)=(t,t^2)$ . Using the above notation, $R(a,b)=b^2-a,
D_1(a,b)=4a+4bi-i$, and $ D_2(a,b)=4a-4bi+i$. By Corollary
\ref{corollary-conics},  for $A\in \cc$ the conchoid is rational. We
analyze the lines given by $D_1$ and $D_2$. We take the
parametrization  ${\cal Q}(h)=(\frac{1}{4}i-h\,i, h)$ of $D_1$. The
rational function $\Delta(h,t)$ is
\[ \Delta(h,t)=\frac{1}{16}(4t^2+4it+1-8h)(2t-i)^2.\]
The discriminant of $4t^2+4it+1-8h$ w.r.t. $t$ is $128(1-4h)$. So,
the only candidate generated by $D_1$ is ${\cal
Q}(1/4)=(0,\frac{1}{4})$ that, indeed, is a double rational focus.
Analyzing $D_2$ one reaches the same point. So, the only double
rational focus for the parabola $\cc$ is  $(0,\frac{1}{4})$ (compare
to Remark \ref{Remark-a-def-rpn} and Example
\ref{example-parabola-rat-focus}); note that we have got the focus
of the parabola.
\end{example}


\begin{example}\label{example-circle-double-rat-focus} {\sf (Double rational focuses for the circle)}
Let $\cc$ be the circle over $\Bbb C$ parametrized by
$$\cP(t)=\left(\frac{2t}{t^2+1},\frac{t^2-1}{t^2+1}\right).$$  One
has that $R(a,b)=a^2+b^2-1$, $D_1(a,b)=-4(a+bi-i)(a+bi)$,
$D_2=-4(a-b i+i)(a-bi)$. Therefore,  $\cal F$ is $\cc$ union the
four lines defined $D_{1,1}:=a+bi,  D_{1,2}:=a+bi-i, D_{2,1}:=a-bi,
D_{2,2}:=a-bi+i$ . By Corollary \ref{corollary-conics},  we only
need to analyze the lines given by $D_{i,j}$. We take the
parametrization  ${\cal Q}(h)=(-i h, h)$ of $D_{1,1}$. Then
\[ \Delta(h,t)=-(-t-i+2th-2ih)(t+i)(t-i)^2.\]
The resultant  of  $(-t-i+2th-2ih)(t+i)$ and its derivative w.r.t.
$t$ is $16h^2(-1+2h))$. So, the  candidates generated by $D_{1,1}$
are ${\cal Q}(0)=(0,0)$ and ${\cal Q}(1/2)=(-\frac{i}{2},
\frac{1}{2})$. One checks that $(0,0)$ is double rational but
$(-\frac{i}{2}, \frac{1}{2})$ is not. Reasoning with the other three
lines no new focuses are found. So, the only double rational focus
is the center of the circle.
\end{example}


\begin{example}\label{example-ellipse-double-rat-focus} {\sf (Double rational focuses for the ellipse)}
Let $\cc$ be the ellipse over $\Bbb C$ parametrized by
$$\cP(t)=\left(\frac{4t}{t^2+1},\frac{3(t^2-1)}{t^2+1}\right).$$
One has that $R(a,b)=a^2/4+b^2/9-1$,
$D_1(a,b)=-4D_{1,1}D_{1,2}D_{1,3}$,  $D_2=-4D_{2,1}D_{2,2}D_{2,3}$,
where $D_{1,1}=a+b i+\sqrt{5}i, D_{1,2}=a+b i-\sqrt{5}i,
D_{1,3}=a+bi-3i$, and $D_{2,j}$ is the conjugate polynomial of
$D_{1,j}$. By Corollary \ref{corollary-conics},  we only need to
analyze the lines given by $D_{i,j}$. We take the parametrization
${\cal Q}(h)=(-i h-i \sqrt{5}, h)$ of $D_{1,1}$. Then
\[ \Delta(h,t)=\frac{- ( \sqrt {5}+3 )}{4}  ( -3\,{t}^{2}-4\,it+3+\sqrt {5
}+2\,h{t}^{2}+2\,h+{t}^{2}\sqrt {5} )  ( 2\,t-3\,i+i\sqrt { 5}
)^{2}.\] The resultant  of  $3 {t}^{2}-4 it+3+\sqrt {5}+2 h{t}^{2}+2
h+{t}^{2}\sqrt {5}$ and its derivative w.r.t. $t$ is
\[32\, ( 5/2-3/2\,h-3/2\,\sqrt {5}+{h}^{2}+3/2\,h\sqrt {5}) h. \]
 So, the  candidates generated by $D_{1,1}$ are (let $\alpha=3/2-1/2\sqrt {5}$)
$${\cal Q}(\alpha)= (-i \alpha -i\sqrt {5},\alpha), \, {\cal Q}(-\sqrt{5})=(0,-\sqrt {5}),\,\,{\cal Q}(0)=(-i\sqrt {5},0).$$ One checks that
${\cal Q}(-\sqrt{5})$ and ${\cal Q}(0)$ are double rational but
${\cal Q}(\alpha)$ is not. Using $D_{1,2}$ one deduces that
$(0,\sqrt {5})$ and $(i\sqrt {5},0)$ are also double rational
focuses. Reasoning with the other  lines no new focuses are found.
So, the only double rational focus are $\{(0,\pm \sqrt{5}), (\pm i
\sqrt{5}) \}$. Note that the focuses of the ellipse have appeared.
\end{example}


\begin{example}\label{example-hyperbola-double-rat-focus} {\sf (Double rational focuses for the hyperbola)}
Let $\cc$ be the hyperbola over $\Bbb C$ parametrized by
$$\cP(t)=\left({\frac {-1-{t}^{2}}{-1+{t}^{2}}},{\frac {2t}{-1+{t}^{2}}}\right).$$
One has that $R(a,b)=a^2-b^2-1, D_1(a,b)=-4D_{1,1}D_{1,2}D_{1,3}$,
$D_2=-4D_{2,1}D_{2,2}D_{2,3}$, where $D_{1,1}=a+b i-\sqrt{2},
D_{1,2}=a+bi +\sqrt{2}, D_{1,3}=a+bi+1$, and $D_{2,j}$ is the
conjugate polynomial of $D_{1,j}$. Reasoning as before, one deduces
that the double rational focuses are $(\pm \sqrt{2},0), (0,\pm i
\sqrt{2})$. Note that Note that the focuses of the hyperbola have
appeared.
\end{example}

\para

\noindent {\sf \underline{Detecting  rational focuses.}} For
analyzing the existence, and actual computation, of rational focuses
we  apply Theorem \ref{theorem-irreducible}. Therefore, we consider
a proper parametrization $\cP(t)$ of $\cc$ and we analyze the
rationality of $\gp$. For this purpose, we analyze the genus of
$\gp$ in terms of the parameters $a,b$ that define the focus. In the
following example, we illustrate these ideas in the case of the
parabola.


\begin{example}\label{example-parabola-rat-focus} {\sf (Rational focuses for the parabola)}
Let $\cc$ be the parabola over $\Bbb C$ parametrized by
$\cP(t)=(t,t^2)$. Let
$$g(x_1,x_2,a,b)= -2\,{ x_2}\,{x_{1}}\, +2\,{x_2}\,a +x_{2}^{2} x_{1}^{2}-x_{2}^{2}b -x_{1}^{2} +b,$$
where $A=(a,b)$ is generic and  let  $G(x_1,x_2,x_3,a,b)$ be the
homogenization of $g$ in the variables $x_1,x_2$. We first observe
that $g$ is primitive w.r.t. $x_2$ iff $A\not\in \cc$. On the other
hand, by Corollary \ref{corollary-conics}, every focus on $\cc$ is
rational. So, we can assume w.l.o.g. that $A\not\in \cc$, and hence
that $g$ is primitive w.r.t. $x_2$. We also know, by Lemma
\ref{lemma-gp} and Example \ref{example-parabola-double-rat-focus},
that the primitive part w.r.t. $x_2$ of $g$ is irreducible iff
$A\neq (0,1/4)$.

\para

We start analyzing the points at infinity.
$G(x_1,x_2,0,a,b)=x_{1}^{2}x_{2}^{2}$. Thus, the points at infinity
of $\gp$ are $P_1:=(1:0:0)$ and $P_2:=(0:1:0)$ independently on $A$.
Moreover $P_1$, $P_2$ are, independently on $A$, double points.
Moreover, $P_1$ is always ordinary (the tangents are given by
$x_{2}^{2}-x_{3}^{2}$) and, if $b\neq 0$, $P_2$ is ordinary too (the
tangents are given by $x_{1}^{2}-b x_{3}^{2}$). Now, we analyze the
affine singular locus. For this purpose, we compute a reduced
Gr\"obner basis $\cal G$ of $\{g, \partial g/\partial x_1, \partial
g/\partial x_2\}$, as polynomials in ${\Bbb C}(a,b)[x_1,x_2]$,
w.r.t. the  graded reverse lexicographic  order with $x_1<x_2$, The
basis is  ${\cal G}=\{1\}$ and
\[ 1=\frac{a_1(x_1,x_2,a,b)}{c(a,b)} g+ \frac{a_2(x_1,x_2,a,b)}{c(a,b)} \frac{\partial g}{\partial x_1}+\frac{a_3(x_1,x_2,a,b)}{c(a,b)} \frac{\partial g}{\partial x_3}, \]
where

\vspace{1 mm}

\noindent $
a_1(x_1,x_2,a,b)=12\,x_{2}\,x_{1}\,{a}^{2}+8\,x_{2}\,x_{1}\,{b}^{2}
-6\,x_{2}\,x_{1}\,b+x_{2}\,x_{1}-8\,x_{2}^{2}b{a}^{2}
+4\,x_{2}^{2}{b}^{2}-x_{2}^{2}b-4\,ax_{1}\,b-ax_{1}+8
\,x_{2}\,{a}^{3}+8\,x_{2}\,a{b}^{2}-2\,x_{2}\,ab+x_{2}\,a+
8\,{a}^{2}b-14\,{a}^{2}-16\,{b}^{2}+8\,b-1, $

\vspace{2 mm}

\noindent $
a_2(x_1,x_2,a,b)=-6\,x_{1}^{2}x_{2}\,{a}^{2}-4\,x_{1}^{2}x_{2}\,{b}^{2}+3\,x_{1}^{2}x_{2}\,b-1/2\,x_{1}^{2}x_{2}+2\,a
x_{1}^{2}b+1/2\,a
x_{1}^{2}+4\,{a}^{3}x_{1}\,x_{2}-2\,ax_{1}\,x_{2}\,b+1/2
\,ax_{1}\,x_{2}-4\,x_{1}\,{a}^{2}b+7\,x_{1}\,{a}^{2}+8\,x_{1}\,{b}^{2}-4\,x_{1}\,b+1/2\,x_{1}+2\,
bx_{2}\,{a}^{2}+4\,{b}^{3}x_{2}-{b}^{2}x_{2}-8\,{a}^{3}-6\,a{b}^{2}+9/2\,ab-1/2\,a,
 $

\vspace{2 mm}

\noindent $  a_3(x_1,x_2,a,b)=4\,b
x_{2}^{3}{a}^{2}-2\,{b}^{2}x_{2}^{3}+1/2\,b
x_{2}^{3}-8\,{a}^{3}x_{2}^{2}-4\,a x_{2}^{2}{b}^{2}+3\,a x_{2}^{2
}b-a x_{2}^{2}-4\,bx_{2}\,{a}^{2}+6\,x_{2}\,{a}^{2}+6\,{b}^{
2}x_{2}-7/2\,bx_{2}+1/2\,x_{2}+8\,{a}^{3}+4\,a{b}^{2}-5\,ab+1 /2\,a,
$

\vspace{2 mm}

\noindent $c(a,b)= ( -b+{a}^{2}  ) (
16\,{b}^{2}-8\,b+1+16\,{a}^{2}).$

\vspace{2 mm}

\noindent Moreover, $g, \partial g/\partial x_1, \partial g/\partial
x_2$ are monic w.r.t. the above order. So, by exercise 7 page 284 in
\cite{Cox},  if $c(a,b)\neq 0$ the Gr\"obner basis specializes
properly. Then, let $\cal W$ be the curve in ${\Bbb C}^2$ defined by
$b ( -b+{a}^{2}  ) ( 16\,{b}^{2}-8\,b+1+16\,{a}^{2}).$ So $\cal W$
is the parabola $\cc$ union the lines $b=0$ (call it $\cal L$), and
$a\pm (1/4-b)i=0$ (call them ${\cal L}^{\pm}$). We distinguish
several cases in our analysis:
\begin{itemize}
\item[(1)] If $A\not\in \cal W$, the genus of $\gp$ is 1; note that $(0,1/4)\in \cal W$, and hence $\gp$ is irreducible. Therefore, $A$ is not rational.
\item[(2)] If $A\in \cc$, we already know that $A$ is rational.
\item[(3)] Let $A\in {\cal L}^{+}$.
If $b=0$, then $A=(-1/4 i, 0)$ and $\gp$ has genus 0. So $(-1/4 i,
0)$ is rational.
 Let $b\neq 0$. Then, $P_1$ and $P_2$ are ordinary double points. So, we only need to analyze the affine singular locus. It holds that, independently on $b$, $(-i/2,i)$ is a double point of $\gp$. So, if $a\neq 0$ (i.e. $A\neq (0,1/4)$), $\gp$ has genus 0. Therefore, every $A\in {\cal L}^{+}\setminus \{(0,1/4)\}$ is rational.
\item[(4)] Let $A\in {\cal L}^{-}$. Reasoning as above, one gets that every $A\in {\cal L}^{-}\setminus \{(0,1/4)\}$ is rational.
 \item[(5)] Let $A\in {\cal L}$.  If $A$ is also on another component of $\cal W$, we already know the classification. So, we assume w.l.o.g. that $A\neq (0,0), A\neq (\pm i/4,0)$. Therefore, in this case, $\gp$ has no affine singularity. Thus, we only need to analyze $P_2$ that now is a non-ordinary singularity. For this purpose, we blow up the curve at $P_2$ (see e.g. Chapter 3 in \cite{libro}).  Applying a suitable projective linear change of coordinates, for instance $\{x_1=x_{1}^{*}-x_{3}^{*},x_2=x_{2}^{*}-x_{3}^{*},x_{3}=x_{1}^{*}+x_{3}^{*}\}$, and the Cremona transformation, one deduces that if $a\neq 0$ then $P_2$ has none neighboring singularities. So, in this case the genus of $\gp$ is 1, and thus no new rational focuses appear.
\end{itemize}
Summarizing, jointly with Example
\ref{example-parabola-double-rat-focus}, one has the following table

\para

\begin{center}
\begin{tabular}{|c|c|c|}
  \hline
   & \mbox{{\sf Double rational focuses}} & \mbox{{\sf Rational focuses}}  \\
\hline
$\begin{array}{c} \mbox{{\sf Parabola}} \\ (t,t^2) \end{array}$ &  $\left(0,\frac{1}{4}\right)$ & $\begin{array}{l} (a,a^2),  \,\,a\in \Bbb C \\
\left(\pm \left(\frac{1}{4}-b\right)i,b)\right), \,\, b\in{\Bbb
C}\setminus \{\frac{1}{4}\} \end{array}$
\\
  \hline
\end{tabular}
\end{center}

\end{example}

\end{document}